\documentclass[a4paper,11pt]{amsart}
\usepackage{fancyhdr}
\usepackage{amsmath}
\usepackage{dsfont}
\usepackage{hyperref}
\usepackage[mathscr]{eucal}
\usepackage[cp1251]{inputenc}
\usepackage[english]{babel}
\usepackage{cite,enumerate,float,indentfirst}
\usepackage{graphicx}
\usepackage{xcolor}
\usepackage{latexsym,a4,mathrsfs,amsthm,amsmath,amssymb,url}
\usepackage{amsfonts}
\usepackage{amssymb}

\numberwithin{equation}{section}
\newtheorem{theorem}{Theorem}[section]
\newtheorem{lemma}[theorem]{Lemma}

\theoremstyle{remark}

\begin{document}

\vspace{1in}

\title[Unitary  Similarity]{\bf Unitary Similarity of Nonderogatory Matrices}

%\vspace{2in}

\author[Yu. Nesterenko]{Yu. Nesterenko}
\address{Yu. Nesterenko\\
 Faculty of Computational Mathematics and Cybernetics\\
 Moscow State University \\
 Moscow, Russia
 }
\email{y\_nesterenko@mail.ru}

\keywords{Canonical forms, unitary similarity}

\date{\today}

\begin{abstract}
This paper is dedicated to the problem of verification of matrices for unitary similarity.
For the case of nonderogatory matrices, we have been able to present the new solution for this problem based on geometric approach.
The main advantage of this approach is stability with respect to errors in the initial upper triangular matrix. Since an upper triangular form is usually obtained by approximate methods (e.g. by QR algorithm), the mentioned advantage seems even more significant and allows us to propose the numerically stable and efficient method for verification of matrices for unitary similarity.
\end{abstract}

\maketitle
\tableofcontents

%\selectlanguage{english}

\newpage

\section{Introduction}

Matrices $ A, B \in \mathds{C}^{n \times n}$ are unitarily similar if a
similarity transformation between them can be implemented
using a unitary matrix $U$:
\begin{equation}
B = UAU^{*}.
\end{equation}

A matrix $A \in \mathds{C}^{n \times n}$ is called nonderogatory if its Jordan blocks have
distinct eigenvalues. Equivalently, a matrix $A \in \mathds{C}^{n \times n}$ is nonderogatory if and only if
its characteristic polynomial and minimum polynomial coincide.

This paper concerns the verification of matrices for
unitary similarity. Based on other authors' works concerning
this problem, two basic approaches can be identified.

In the first, a complete system of matrix invariants
under a unitary similarity transformation is constructed.
In a sense, the final result in this direction is
the Specht-Pearcy criterion (see~\cite{Specht1, Pearcy1}), which
reduces the question to verifying conditions of the form
\begin{equation}
\textrm{tr}\, W(A,A^*) = \textrm{tr}\, W(B,B^*)
\end{equation}
for all words $W(s, t)$ of length at most $2n^2$. However, it
seems that the number of words to be verified is
strongly overestimated (see~\cite{Mumagh1, Sibir1, Laffey1, Bhatta1}). Moreover, this
method cannot find a matrix generating a given unitary similarity.

The second approach is free of this shortcoming
and consists of constructing a canonical form of matrices
with respect to unitary similarity transformations.
Inductive definitions of the canonical form of a matrix
were proposed in~\cite{Brenner1, Little1, Radjavi1},
but it is hard to visualize the final canonical form.
In more recently work~\cite{Futorny1} the autors, considered the set of nonderogatory matrix,
constructed more visual canonical form.

In this work we poropose the geometric approach to solving the problem for nonderogatory matrices.
Given an arbitrary nonderogatory matrix, we construct a finite family
of unitarily similar matrices for it (this family is called
canonical). Whether or not two matrices are unitarily
similar can be answered by verifying the intersection of
their corresponding families.
This method for unitary similarity verification has the significant advantage
over the method~\cite{Futorny1} based on construction of the canonical form,
it is stable with respect to errors in the initial matrix.
This last aspect is discussed at the end of this paper.

While constructing a canonical family, we start
from an upper triangular matrix form. Specifically, by
the Schur theorem, any matrix $A \in \mathds{C}^{n \times n}$ can be reduced to such
a form by using a unitary similarity transformation:
\begin{equation}
\Delta = \begin{bmatrix}
\lambda_1 & \Delta_{12}  & \Delta_{13} & \ldots      & \Delta_{1n} \\
\quad     & \lambda_2    & \Delta_{23} & \ldots      & \Delta_{2n} \\
\quad     & \quad        & \lambda_3   & \quad       & \quad       \\
\quad     & \quad        & \quad       & \ddots      & \quad       \\
\quad     & \quad        & \quad       & \quad       & \lambda_n
\end{bmatrix},
\end{equation}
where $\lambda_1,\ldots,\lambda_n$ are eigenvalues of the matrix with multiplicity in some fixed order.
The next statement let us to restrict the set of unitary transformations while operating with a nonderogatory triangular matrices.

\begin{lemma}
Let $A$ be a nonderogatory complex $n \times n$ matrix, and let $\Delta$ be its upper triangular form with
eigenvalues $\lambda_1,\ldots,\lambda_n$ on the diagonal in some fixed order. Then the magnitudes of the elements $\Delta_{ij}$, where $i < j$,
are uniquely determined.
\end{lemma}

\textbf{Proof.} Since the similar proposition for matrices with simple eigenvalues is known~\cite{Little1}, we can consider the case of nonderogatory matrix
with the single eigenvalue $\lambda$.
Let $\Delta$ is obtained from $A$ by unitary similarity transformation
\begin{equation}\label{eq:2}
\Delta = Q^{*}AQ
\end{equation}
where $Q = (q_{1} q_{2} \ldots q_{n})$ is unitary matrix.
Rewriting equation \eqref{eq:2} as $AQ = Q\Delta$, one may see that
$q_{1}$ is normalized eigenvector of the matrix $A$ corresponding to the eigenvalue $\lambda$.
Further,
\begin{equation}\label{eq:3}
Aq_{2} = \Delta_{12}q_{1} + \lambda q_{2}
\end{equation},
hence $(A - \lambda E)q_{2} = 0$, i.e. $q_{2}$ is generalized eigenvector of A. Adding the condition of orthonormality of the pair $q_{1}, q_{2}$,
one obtains that $q_{2}$ is uniquely determined up to multiplication by a scalar of unit modulus.
Continuing in the same vein, we can see that the matrix $Q$ is uniquely determined up to multiplication by a diagonal unitary matrix,
but such a transformations preserve the magnitudes of the off-diagonal elements of upper triangular form $Q$.
Thus the lemma is proved.

Using the last lemma we can limit our consideration to studying the action of the group of unitary similarity transformations
with diagonal matrices on the set of upper triangular matrices:
\begin{equation}
\Delta \mapsto X\Delta X^{*}, \quad X = diag(e^{i \psi_1},\ldots,e^{i \psi_{n-1}},1),
\end{equation}
\begin{equation}
\{X \Delta X^*\}_{ij} = \begin{cases}
\Delta_{ij} e^{i (\psi_i - \psi_i)}  & \quad i < j < n,\\
\Delta_{ij} e^{i \psi_i}  & \quad i < j = n,\\
\lambda_{i} & \quad i = j,\\
0 & \quad i > j
\end{cases}
\end{equation}
(assumming the last diagonal entry of $X$ is $1$, we remove a scalar factor from $X$).

\section{Preliminary constructions}

Let $M$ denote the range of the parameters of the
matrix
\begin{equation}
M = \{(r_{12},\ldots,r_{n-1,n};\varphi_{12},\ldots,\varphi_{n-1,n}),
\quad r_{ij},\varphi_{ij} \in \mathbb{R} \},
\end{equation}
and let $M_{r}$ denote its restriction for fixed $r_{ij}$:
\begin{equation}
M_{r} = \{(r;\varphi)\in M : r_{ij} \text{are fixed} \}.
\end{equation}

The indices $i$ and $j$ run over the values $1 \leq i < j \leq n$
and are ordered lexicographically. For elements of $M$
and $M_{r}$, several equivalent forms of notation are used:
\begin{equation}
(r_{12},\ldots,r_{n-1,n};\varphi_{12},\ldots,\varphi_{n-1,n}) =
(r;\varphi_{12},\ldots,\varphi_{n-1,n}) = (r;\varphi).
\end{equation}

Looking ahead, $r_{ij}$ and $\varphi_{ij}$ will later play the role of
absolute values and arguments of off-diagonal elements of $\Delta$.
Despite this geometric interpretation, no constraints
are as yet imposed on $r_{ij}$ and $\varphi_{ij}$ and the indetermination of
$\varphi_{ij}$ at $r_{ij} = 0$ is ignored. At this stage, we work with the
formally defined range $M$.

On $M$ we introduce the family of transformations
\begin{equation}
X_{\psi}: (r; \varphi) \mapsto (r; \tilde{\varphi}),
\end{equation}
\begin{equation}
\begin{split}
&\tilde{\varphi}_{ij} = \varphi_{ij} + \psi_i - \psi_j, \quad 1 \leq i < j \leq n-1,\\
&\tilde{\varphi}_{in} = \varphi_{in} + \psi_i, \quad 1 \leq i \leq
n-1.
\end{split}
\end{equation}
Each such a transformation is defined by a parameter
vector $\psi = (\psi_{1},\ldots,\psi_{n-1}) \in \mathbb{R}^{n-1}$.

Consider a subset of matrices $K \subset M$ whose elements satisfy the system of equations
\begin{equation}\label{eq:33}
-\sum_{k=1}^{s-1} r_{ks}\varphi_{ks} + \sum_{k=s+1}^{n}
r_{sk}\varphi_{sk} = 0, \quad s = 1,\ldots,n-1.
\end{equation}

The summation indices in \eqref{eq:33} are visually described
by the diagram
\begin{equation}
\begin{bmatrix}
\lambda_1             & \quad     & \quad   & *          & \quad          & \quad   \\
\quad                 & \ddots    & \quad   & *          & \quad          & \quad    \\
\quad                 & \quad     & \ddots  & *          & \quad          & \quad     \\
\quad                 & \quad     & \quad   & \lambda_s  & *              & *          \\
\quad                 & \quad     & \quad   & \quad      & \ddots         & \quad       \\
\quad                 & \quad     & \quad   & \quad      & \quad          & \lambda_n    \\
\end{bmatrix}.
\end{equation}

The reduction of an arbitrary matrix of $M$ to a $K$
form by applying a transformation $X_{\psi}$ is reduced to
finding the parameters of this transformation $\psi = (\psi_1,\ldots,\psi_{n-1})$
by solving the system of linear equations
\begin{equation}\label{eq:34}
R(r) \psi = - b(r,\varphi)
\end{equation}
with the symmetric matrix
\begin{equation}
\{R(r)\}_{ij} = \begin{cases}
-r_{ij} & \quad i < j,\\
\sum_{k=1}^{i-1} r_{ki} + \sum_{k=i+1}^{n} r_{ik} & \quad i = j,\\
-r_{ji} & \quad i > j
\end{cases}
\end{equation}
and with a righthand side that is linear in $r$ and $\varphi$:
\begin{equation}
\begin{split}
&b(r,\varphi) = (b_{1}(r,\varphi),\ldots, b_{n-1}(r,\varphi)) ,\\
&b_{s}(r,\varphi) = -\sum_{k=1}^{s-1} r_{ks}\varphi_{ks} +
\sum_{k=s+1}^{n} r_{sk}\varphi_{sk}, \quad s = 1,\ldots,n-1.
\end{split}
\end{equation}

System \eqref{eq:34} has some remarkable properties.

\begin{theorem} \label{th:1}
(i) For any $\varphi_{ij}$ and nonnegative $r_{ij}$,
system \eqref{eq:34} has a solution; i.e.,
\begin{equation}
- b(r,\varphi) \in \textrm{Im} R(r), \quad \forall \varphi_{ij},
\quad \forall r_{ij} \geq 0.
\end{equation}

(ii) For all $r_{ij} \geq 0$, the determinant $det R(r) \neq 0$ is
nonzero if and only if the indices of the nonzero elements
$r_{ij} > 0$ contain a collection $(ij)^{1},\ldots,(ij)^{n-1}$ such that the
set $\{~\psi_{i^p}~-~\psi_{j^p} \quad (\text{respectively}
\quad \psi_{i^p}, \quad \text{if} \quad j^p = 0), \quad p =
1,\ldots,{n-1} \}$ forms a linearly independent system of functions of
variables $(\psi_{1},\ldots,\psi_{n-1})$.

(iii) Even if $det R(r) = 0$ for some $r_{ij} \geq 0$, the solution
of the equation $\psi = (\psi_1,\ldots,\psi_{n-1})$ is such that the
quantities
\begin{equation}
\begin{split}
&r_{ij}(\psi_i - \psi_j), \quad 1 \leq i < j \leq n-1 \text{\quad and}\\
&r_{in}\psi_i, \quad 1 \leq i \leq n-1
\end{split}
\end{equation}
are uniquely defined. This means that nonuniqueness in
the definition of $\psi_i - \psi_j$ occurs if and only if $r_{ij} = 0$.
\end{theorem}

\textbf{Proof.} On the set $M_{r}$, we introduce the natural
structure of a Euclidean space:
\begin{equation}
(r;\varphi^{(1)}) + (r;\varphi^{(2)}) = (r;\varphi^{(1)}_{12} +
\varphi^{(2)}_{12},\ldots,\varphi^{(1)}_{n-1,n} +
\varphi^{(2)}_{n-1,n}),
\end{equation}
\begin{equation}
\alpha (r;\varphi) = (r;\alpha \varphi_{12},\ldots,\alpha
\varphi_{n-1,n}), \quad \alpha \in \mathbb{R},
\end{equation}
\begin{equation}
\langle(r;\varphi^{(1)}),(r;\varphi^{(2)})\rangle = \sum_{1 \leq i <
j \leq n} \varphi^{(1)}_{ij} \varphi^{(2)}_{ij}.
\end{equation}

Then $K_{r} = M_{r} \bigcap K$ is a linear subspace of $M_{r}$ that is
orthogonal to all linear manifolds of the form
\begin{equation}
G_{r,\varphi} = \{ (r;\varphi) + \sum_{1 \leq i < j \leq n-1}
(\psi_{i} - \psi_{j}) \, r_{ij} \, \mathbb{I}_{ij} + \sum_{1 \leq i
\leq n-1} \psi_{i} \, r_{in} \, \mathbb{I}_{in}, \quad \psi \in
\mathbb{R}^{n-1} \},
\end{equation}
\begin{equation}
\mathbb{I}_{ij} = (r;0,\ldots,0,\stackrel{(ij)}{1},0,\ldots,0).
\end{equation}

The dimensions of $K_{r}$ and $G_{r,\varphi}$ depend on $r$, but
their sum is a constant:
\begin{equation}
\textrm{dim} \, K_{r} + \textrm{dim} \, G_{r,\varphi} = \textrm{dim}
\, M_{r}.
\end{equation}
In other words, in the Euclidean space $M_{r}$, the linear
space $K_{r}$ and the linear manifold $G_{r,\varphi}$ are mutually
orthogonal and the sum of their dimensions is the
complete one. This implies that they have a unique
intersection point $(r;\varphi') = K_{r} \bigcap G_{r,\varphi}$. This intersection
condition corresponds to the system of equations
\begin{equation}\label{eq:35}
R(r^2) \psi = - b(r,\varphi),
\end{equation}
where $r^2$ denotes the vector
\begin{equation}
r^2 = (r^2_{12},\ldots,r^2_{n-1,n}).
\end{equation}
In terms of $\psi$, the existence and uniqueness of an
intersection point $(r;\varphi')$ means that system \eqref{eq:35} is
solvable with arbitrary $\varphi_{ij}$ and $r_{ij}$ and that the values
\begin{equation}\label{eq:36}
\begin{split}
&r_{ij}(\psi_i - \psi_j) = f_{ij}, \quad 1 \leq i < j \leq n-1 \text{\quad and}\\
&r_{in}\psi_i = f_{in}, \quad 1 \leq i \leq n-1
\end{split}
\end{equation}
are uniquely determined from it.

Assume that there exists an index set $(ij)^{1},\ldots,(ij)^{n-1}$
corresponding to the nonzero elements of $R(r^2)$ such
that the set $\{~\psi_{i^p}~-~\psi_{j^p} \quad
(\text{respectively} \quad \psi_{i^p}, \quad \text{if} \quad j^p
= 0), \quad p = 1,\ldots,n-1 \}$ forms a linearly independent system of
functions of variables $(\psi_{1},\ldots,\psi_{n-1})$.
Then a nondegenerate system of linear equations can be composed
of relations \eqref{eq:36} and $\psi = (\psi_1,\ldots,\psi_{n-1})$ can
be uniquely determined. Thus, under the conditions
formulated, system \eqref{eq:35} has a unique solution and,
hence, $det R(r) \neq 0$. The converse can be proved by
contradiction.

The above results are extended to system \eqref{eq:34} by
making the substitution $r'_{ij} = r^2_{ij} \geq 0$. The proof is complete.

Returning to the matrix $\Delta$, we use Theorem~\ref{th:1}
to construct the family of matrices that are unitarily similar to $\Delta$.

With the help of the elements of $\Delta$, we set up the
system of linear equations
\begin{equation}\label{eq:37}
R(r) \psi = - b(r,\varphi + 2\pi m),
\end{equation}
where $r$, $\varphi$, and $m$ are defined as
\begin{equation}\label{eq:38}
r_{ij} = |\Delta_{ij}|, \quad \varphi_{ij} = \textrm{arg} \, \Delta_{ij} -
\pi, \quad m_{ij} \in \mathbb{Z}
\end{equation}
Note that, despite the indetermination of $\varphi_{ij}$ at $r_{ij} = 0$,
the system of equations is uniquely defined.

Solving this system for $\psi = (\psi_1,\ldots,\psi_{n-1})$, we
construct the matrix $X \Delta X^*$, $X = diag(e^{\textrm{i}
\psi_1},\ldots,e^{\textrm{i} \psi_{n-1}},1)$,
which is unitarily similar to the original one. Again, if
for some matrix $\Delta$ the parameter vector $\psi$ is not
determined uniquely from system \eqref{eq:37}, then, by Theorem~\ref{th:1},
this nonuniqueness is such that the matrix
$X \Delta X^{*}$ is uniquely determined.

The matrix generated by this procedure from $\Delta$
with the parameter vector $m = (m_{12},\ldots,m_{n-1,n})$ is
denoted by $\mathcal{K}(\Delta,m)$.

\section{Algorithm for constructing the canonical family}

Now, we consider two nonderogatory upper triangular matrices
$\Delta^{(1)}$ and $\Delta^{(2)}$ with identical sets of eigenvalues. The
eigenvalues are assumed to be identically ordered on the matrix diagonals. For these
matrices, we introduce $r^{(1)}_{ij}$, $\varphi^{(1)}_{ij}$ and
$r^{(2)}_{ij}$, $\varphi^{(2)}_{ij}$ similar
to \eqref{eq:38}. The matrices $\Delta^{(1)}$ and $\Delta^{(2)}$ are
related by a unitary similarity transformation if and only if

(i) $r^{(1)}_{ij} = r^{(2)}_{ij}, \quad 1 \leq i < j \leq n$ and

(ii) there exist sets $(\psi_{1},\ldots,\psi_{n-1}) \in
\mathbb{R}^{n-1}$ and $(k_{12},\ldots,k_{n-1,n}) \in
\mathbb{Z}^{\frac{n(n-1)}{2}}$ such that, for indices $(ij)$
corresponding $r^{(1)}_{ij} = r^{(2)}_{ij} > 0$, we have
\begin{equation}\label{eq:39}
\varphi^{(1)}_{ij} + 2\pi k^{(1)}_{ij} = \varphi^{(2)}_{ij} + 2\pi
k^{(2)}_{ij} + \psi_{i} - \psi_j.
\end{equation}

This implies that a unitary similarity of $\Delta^{(1)}$ and $\Delta^{(2)}$
is equivalent to $\mathcal{K}(C^{(1)},k^{(1)}) = \mathcal{K}(C^{(2)},k^{(2)})$ for some
integer parameter vectors $k^{(1)} = (k^{(1)}_{12},\ldots,k^{(1)}_{n-1,n})$ and
$k^{(2)} = (k^{(2)}_{12},\ldots,k^{(2)}_{n-1,n})$.

Let us represent the above criterion in an effective
form. Define a subset $\mathbb{I} \subset \mathbb{Z}^{\frac{n(n-1)}{2}}$:
\begin{multline}
\mathbb{I} = \{ k \in \mathbb{Z}^{\frac{n(n-1)}{2}}: \quad k_{ij} =
0, \pm 1, \quad 1 \leq i < j \leq n-1, \\ k_{in} = 0, \quad 1 \leq i
\leq n-1 \}.
\end{multline}

\begin{theorem} \label{th:2}
The matrices $\Delta^{(1)}$ and $\Delta^{(2)}$ are unitarily similar
if and only if there exist vectors $k^{(1)}, k^{(2)} \in \mathbb{I}$ such that $\mathcal{K}(C^{(1)},k^{(1)}) = \mathcal{K}(C^{(2)},k^{(2)})$.
\end{theorem}

\textbf{Proof.} Let $\Delta^{(1)}$ and $\Delta^{(2)}$ be unitarily similar and all
their elements above the diagonal be nonzero. Then, as was
shown above, there exist vectors $(\psi_{1},\ldots,\psi_{n-1}) \in \mathbb{R}^{n-1}$ and
$k^{(1)}, k^{(2)} \in \mathbb{Z}^{\frac{n(n-1)}{2}}$ such that equalities \eqref{eq:39}
hold for all $(ij)$. We use them to make up the following linear combinations:
\begin{multline}
\varphi^{(1)}_{ij} - \varphi^{(1)}_{in} + \varphi^{(1)}_{jn} + 2\pi
(k^{(1)}_{ij} - k^{(1)}_{in} + k^{(1)}_{jn}) = \\ =
\varphi^{(2)}_{ij} - \varphi^{(2)}_{in} + \varphi^{(2)}_{jn} + 2\pi
(k^{(2)}_{ij} - k^{(2)}_{in} + k^{(2)}_{jn}).
\end{multline}
One may see that the $\psi$-dependent terms have canceled out.
A feature of these linear combinations is that
they are invariant under the action of transformations $X_{\psi}$
on the linear space of vectors $\varphi = (\varphi_{12},\ldots,\varphi_{n-1,n})$.
Moreover, these combinations form a basis in the subspace of linear functionals invariant under $X_{\psi}$.

Note that the conditions $\varphi^{(s)}_{ij} \in [-\pi, \pi)$ imply $\varphi^{(s)}_{ij} - \varphi^{(s)}_{in} +
\varphi^{(s)}_{jn} \in (-3\pi, 3\pi)$, which in turn imply the following constraints on $k^{(1)}$ and $k^{(2)}$:
\begin{equation}\label{eq:42}
(k^{(1)}_{ij} - k^{(1)}_{in} + k^{(1)}_{jn}) - (k^{(2)}_{ij} -
k^{(2)}_{in} + k^{(2)}_{jn}) = 0, \pm 1, \pm 2.
\end{equation}

At the same time, the algorithm for deriving the
matrix $\mathcal{K}(\Delta,0)$ shows that the arguments of its
elements are linearly expressed in terms
of $\varphi_{ij}$:
\begin{equation}
\tilde{\varphi}_{ij} \in \mathcal{L}
(\varphi_{12},\ldots,\varphi_{n-1,n}),
\end{equation}
Moreover, these linear combinations must be invariant
under $X_{\psi}$, so their form can be refined:
\begin{equation}
\tilde{\varphi}_{ij} \in \mathcal{L} ( \{\varphi_{ij} - \varphi_{in}
+ \varphi_{jn}\}, \quad 1 \leq i < j \leq n-1 ).
\end{equation}
Combining this with \eqref{eq:42}, we obtain the sufficiency of
verifying the equalities $\mathcal{K}(\Delta^{(1)},k^{(1)}) = \mathcal{K}(\Delta^{(2)},k^{(2)})$ for
$k^{(1)}, k^{(2)} \in \mathbb{I}$.

In the presence of zero elements above the diagonal of
$\Delta^{(1)}$ and $\Delta^{(2)}$, the proposition is proved with slight modifications.

The finite set of matrices $\mathcal{K}(\Delta,k)$, $k \in \mathbb{I}$, that are
unitarily similar to $\Delta$ is called the canonical family of
the given matrix.

Thus, the following algorithm is proposed for verifying unitary similarity between nonderogatory matrices A and B with
the same set of eigenvalues:

(i) Reduce these matrices to an upper triangular
form with identically ordered eigenvalues on the diagonal to obtain matrices $\Delta^{(1)}$ and $\Delta^{(1)}$:
\begin{equation}
\Delta^{(1)} = U_1 A U_1^{*}, \quad \Delta^{(1)} = U_2 B U_2^{*}
\end{equation}

(ii) For $\Delta^{(1)}$ and $\Delta^{(1)}$, construct their canonical families
$\mathcal{K}(\Delta^{(1)},k^{(1)})$ and $\mathcal{K}(\Delta^{(2)},k^{(2)})$, $k^{(1)}, k^{(2)} \in \mathbb{I}$.

(iii) If these families intersect for some $k^{(1)}, k^{(2)} \in \mathbb{I}$ and
\begin{equation}
\mathcal{K}(\Delta^{(1)},k^{(1)}) = X_1 \Delta^{(1)} X_1^{*}, \quad
\mathcal{K}(\Delta^{(2)},k^{(2)}) = X_2 \Delta^{(2)} X_2^{*},
\end{equation}
then the original matrices are similar and
\begin{equation}
B = UAU^{*}, \quad U = U_2^* X_2^* X_1 U_1.
\end{equation}
Otherwise, they are not similar.

\section{Numerical stability}

The approach presented above significanly differs from earlier approaches to the problem studied.
As a rule, different approaches (e.g. \cite{Little1, Futorny1}), based on the Schur upper triangular form,
tried to create as many positive elements above the diagonal as possible.
But such a property of a desired canonical form inevitably leads to
the form unstable with respect to errors in initial triangular form.
One may observe the present effect on the next example:
\begin{equation}
A(\varepsilon) = \begin{bmatrix}
1             & i         & i       & i            \\
0             & 2         & i       & i            \\
0             & 0         & 3       & \varepsilon  \\
0             & 0         & 0       & 4
\end{bmatrix},
\end{equation}
where $\varepsilon$ is a complex number. If the initially "strategy" of obtaining
the greatest possible number of positive off-diagonal elements is to start with superdiagonal elements,
then one can chose a $A(\varepsilon)$ arbitrary close (e.g. with respect to the Frobenius norm) to $A(0)$,
but their canonical forms won't satisfy this property.
The stability property seems even more significant due to the fact that usually
an upper trianglular form of a matrix is obtained by approximate methods (e.g. QR algorithm).

From the geometric point of view the constructed canonical family is the finite set of the
ruled surfaces, such that an orbit of each nonderogatory matrix intersects each of them in a single point.
The stability of this set of intersection points follows from the continuity of quantities \eqref{eq:36} determined from system \eqref{eq:37}.
The present property is of special interest in
the context of the result obtained in~\cite{Paulsen1}.
Many ideas used by the author were taken from~\cite{Arnold1}.
Specifically, a minimal continuous extension of a
canonical Jordan form was constructed in~\cite{Arnold1}.
Some of the results presented above are reflected in ~\cite{Nester1}.

I am deeply grateful to Professor Kh.D. Ikramov
for his interest in this work and helpful discussions.

\newpage

\bibliographystyle{plain}
\bibliography{literatur}

\end{document}